\documentclass[12pt]{amsart}
\usepackage{color,graphicx,amsmath,amsfonts,amssymb,amsxtra,setspace,xspace,subfig}
\usepackage[colorlinks=true]{hyperref}
\hypersetup{urlcolor=blue, citecolor=red, linkcolor=blue}
\newtheorem{theorem}{Theorem}[section]
\newtheorem{lemma}[theorem]{Lemma}
\newtheorem{definition}[theorem]{Definition}

\newtheorem{remark}[theorem]{Remark}

\newtheorem{problem}[theorem]{Problem}

\numberwithin{equation}{section}

\begin{document}

\title[Meixner random variables and their quantum operators]
{Meixner random variables and their quantum operators}

\author{Nobuaki Obata}
\address{Center for Data-driven Science and Artificial Intelligence \\
Tohoku University\\
Sendai, 980-8576, Japan\\
E:mail: obata@tohoku.ac.jp}

\author{Aurel I. Stan}
\address{Department of Mathematics \\
Ohio State University at Marion \\
1465 Mount Vernon Avenue \\
Marion, OH 43302, U.S.A.\\
E-mail: stan.7@osu.edu}

\author{Hiroaki Yoshida}
\address{Department of Information Sciences\\
Ochanomizu University\\
2-1-1, Otsuka, Bunkyo\\
Tokyo, 112-8610, Japan\\
E-mail: yoshida@is.ocha.ac.jp}

\thanks{NO is supported in part by
JSPS Grant-in-Aid for Scientific Research No.~23K03126.}

%\keywords{radial solutions, singular potential}
%\subjclass{Primary: 35J25, 35J70, Secondary: 35J60}

\date{\today}

\maketitle

\begin{abstract}
We find the position-momentum decomposition of the quantum operators
of the classic Meixner random variables. The position-momentum
decomposition involves translation operators, which
are used to give a new characterization of the Meixner random variables.
%Finally, the number and second quantization operators, of a constant times the identity
%operators, for the Meixner distributions are described in terms of the translation
%and dilation operators.
\end{abstract}



\section{Introduction}

The field of orthogonal polynomials is a rich and
classic area of investigation that have occupied the minds of many
mathematicians in the last two centuries. 
People have employed various tools like continuous fraction expansions,
generating functions, re-normalization techniques, 
and combinatorial methods to understand the orthogonal polynomials
generated by a random variable (or equivalently, 
by its probability distribution),
having finite moments of all orders.
These investigations focused on the concrete form of
the orthogonal polynomials for various classic random variables.
An equivalent way of studying the underlying orthogonal structure
of the unital algebra generated by a random variable, 
having finite moments of all orders,
that is more operatorial in nature,
and somehow follows the philosophy of quantum mechanics,
has been developed in the last decades by analyzing
the quantum operators: creation, preservation, and annihilation 
operators. 
A particular class of random variables that is well suited
to this operatorial treatment consists of the classic
Meixner random variables which include the shifted and re-scaled 
Gaussian, Poisson, gamma, negative binomial, two parameter hyperbolic,
and binomial distributions.
This class is well suited to this treatment due to the
Lie algebra structure enjoyed by the quantum operators,
that means the fact that the commutator of any two quantum operators
can be expressed as a linear combination of the quantum,
number, and identity operators.
There are two fundamental operators: 
the multiplication by $X$, also called the {\em position} operator,
and the differentiation operator $D = d/dx$, 
also called the {\em momentum} operator. 
Every quantum operator and each linear operator
that maps the space $F$, of all random variables of the form $P(X)$,
where $P$ is a polynomial and $X$ is a given random variable, 
having finite moments of all orders, into itself, can be expressed
uniquely as sum (series) of terms of the form $A_n(X)D^n$,
that means, compositions of powers of the position operator and
powers of the momentum operator. 
In these compositions, the position ``factors" are always placed
at the left of the momentum ``factors". 
The importance of the decomposition into 
a sum (series) of compositions of position and momentum operators
relies on the Heisenberg commutation relation:
\begin{equation}
\left[D, X\right] = I,
\end{equation}
where $I$ is the identity operator.

These facts can be extended to the multi-dimensional case,
where one random variable $X$ is replaced 
by a finite family of random variables
$X_1,X_2,\dots,X_d$, 
defined on the same probability space, 
and having finite moments of all orders. 
On the surface, the random variables $X_1,X_2,\dots,X_d$ commute,
but beneath the surface there are hidden the quantum operators:
creation, preservation, and annihilation,
which are defined in the terms of the orthogonal spaces and 
compositions of the multiplication operators
by $X_1,X_2,\dots,X_d$ (called the position operators) 
and projections onto these orthogonal spaces.
These quantum operators do not commute, 
and the non-zero commutators 
between any two of these quantum operators 
make the study of quantum probability interesting.
In classical probability
a family of random variables is determined 
by their joint probability distribution.
But in quantum probability a family of random variables
is characterized by the concrete form of their joint quantum operators,
the number operator and the second quantization operators of constants
times the identity operator, 
which ultimately lead to a concrete form of
the projection operators on the orthogonal spaces 
which are generated by the family of random variables.

This paper is dedicated to the study of 
the classical Meixner random variables 
along with quantum probability.
As we will see, their quantum operators will be expressed
in terms of translation operators.
Moreover, we will derive a characterization of these random variables 
in terms of their quantum operators.
We will not address the problem of knowing the second quantization
operator of a constant times the identity operator in this paper.
Moreover, we will derive a characterization of these random variables
in terms of these quantum operators.

The paper is structured as follows. 
In section 1 we present a brief background of quantum,
semi-quantum and number operators.
In section 2 we focus on Meixner random variables, 
and find the concrete form of the operators described in Section 1.
In section 3 we characterize the Meixner random variables
using the forms found in Section 2.

\section{Background}

Let $X_1,X_2,\dots, X_d$ be $d\ge1$ real valued random variables,
defined on a common probability space $(\Omega, \mathcal{F},P)$.
We assume that these random variables have finite moments 
of all orders, that means, for all $i \in \{1,2,\dots,d\}$ 
and all $p>0$, we have:
\begin{equation}
E\left[\left|X_i\right|^p\right]
=\int_{\Omega}\left|X_i(\omega)\right|^pdP(\omega)
<\infty.
\end{equation}
Let $F$ be the space of all random variables of the form 
$f(X_1, X_2, \dots, X_d)$, 
where $f$ is a polynomial of $d$ variables with complex coefficients. 
For all $n \in \mathbb{N}\cup\{0\}$
let $F_n$ be the subspace of $F$ consisting of 
all random variables of the form $f(X_1,X_2,\dots, X_d)$, 
where $f$ is a polynomial of degree at most $n$. 
Since the random variables $X_1,X_2, \dots, X_d$ have finite
moments of all orders, we have:
\begin{equation}
\mathbb{C}\cong F_0 \subseteq F_1 \subseteq F_2 
\subseteq \dotsb \subseteq F \subseteq L^2(\Omega,\mathcal{F}, P).
\end{equation}
Because of finite dimensionality $F_n$ is a closed subspace of 
$L^2(\Omega,\mathcal{F},P)$ for all $n \in \mathbb{N} \cup \{0\}$.
For convenience we set $F_{-1}= \{0\}$ (the null space).

For $n \in \mathbb{N} \cup \{0\}$ we define
\begin{equation}
G_n=F_n \ominus F_{n-1},
\end{equation}
that means, $G_n$ is the orthogonal complement of $F_{n-1}$ in $F_n$.
Then obviously,
\begin{equation}\label{02eqn:Fn is sum of Gs}
F_n= G_0 \oplus G_1 \oplus \cdots \oplus G_n.
\end{equation}
We call $G_n$ the {\em $n$th homogenous chaos space} generated by 
$X_1,X_2,\dots,X_d$, and a random variable $f(X_1,X_2, \dots,X_d)$
in $G_n$ a {\em homogenous polynomial random variable} of degree $n$
(here the word ``homogenous" does not have the classic meaning
that all terms have the same degree).

It is essential for us to regard 
$X_1, X_2, \dots, X_d$ not as random variables,
but as the multiplication operators 
\begin{equation}
X_i:f(X_1,X_2,\dots,X_d) \mapsto X_i f(X_1,X_2,\dots,X_d).
\end{equation}
Thus, $X_i$ becomes a linear operator from $F$ into itself.

\begin{lemma}[{\cite[Theorem 1]{an02}, \cite[Lemma 2.1]{aks04}}]
\label{lemma_orthogonal}
Let $i \in \{1,2,\dots,d\}$
and $n \in \mathbb{N}\cup\{0\}$.
The image of $G_n$ by the multiplication operator $X_i$ is
orthogonal to $G_k$ for all $k \neq n-1, n, n+1$,
that is, $X_iG_n \perp G_k$.
\end{lemma}
Since $X_iF_n\subseteq F_{n+1}$ and $G_n \subseteq F_n$,
we see from \eqref{02eqn:Fn is sum of Gs} that
$X_iG_n \subseteq G_0 \oplus G_1 \oplus \cdots \oplus G_{n + 1}$.
Then with the help of Lemma \ref{lemma_orthogonal} we obtain
\begin{equation}
X_iG_n \subseteq G_{n - 1} \oplus G_n \oplus G_{n + 1}.
\end{equation}
Thus, for each homogenous polynomial random variable of degree $n$,
$f \in G_n$, there exist unique $f_{i, n-1} \in G_{n-1}$,
$f_{i,n} \in G_n$, and $f_{i,n+1} \in G_{n+1}$
such that:
\begin{equation}\label{02erqn:orthogonal decomposition of Xif}
X_if = f_{i, n-1} \oplus f_{i,n} \oplus f_{i,n+1},
\end{equation}
where $\oplus$ stands for the orthogonal sum.

According to \eqref{02erqn:orthogonal decomposition of Xif}
we define three operators $D_n^-(i),D_n^0(i)$ and $D_n^+(i)$ by
\begin{equation}
D_n^-(i)f = f_{i,n-1},
\qquad
D_n^0(i)f = f_{i,n},
\qquad
D_n^+(i)f = f_{i, n+1}.
\end{equation}
Then we have
\begin{equation}
D_n^-(i) : G_n \to G_{n-1},
\qquad
D_n^0(i) : G_n \to G_n,
\qquad
D_n^+(i) : G_n \to G_{n+1}.
\end{equation}
The operator $D_n^-(i)$ decreases the degree of
a homogenous polynomial random variable by $1$ unit
and is called an {\em annihilation operator}.
The operator $D_n^0(i)$ does not change the degree of
a homogenous polynomial random variable
and is called a {\em preservation operator}.
Finally, $D_n^+(i)$ increases the degree of
a homogenous polynomial random variable by $1$ unit
and is called a {\em creation} operator.
With these notations Lemma \ref{lemma_orthogonal} can be written as

\begin{lemma}\label{lemma_restricted_quantum_decomposition}
For all $i \in \{1,2,\dots,d\}$ and $n \in\mathbb{N}\cup\{0\}$ 
we have
\begin{equation}
X_i|_{G_n} = D_n^-(i) + D_n^0(i) + D_n^+(i),
\end{equation}
where $X_i|_{G_n}$ stands for the restriction of the 
multiplication operator $X_i$ to $G_n$.
\end{lemma}

We extend the definition of the creation, preservation, 
and annihilation operators to the space $F$ of 
all polynomial random variables in the following manner.
For all $f \in F$ there exist unique $f_0 \in G_0$, 
$f_1 \in G_1$, $f_2 \in G_2$, $\cdots$, with only finitely many 
of them being not zero, such that
\begin{equation}
f = f_0 \oplus f_1 \oplus f_2 \oplus \cdots.
\end{equation}
We define the {\em creation operator} $a^+(i)$ by
\begin{equation}
a^+(i)f = D_0^+(i)f_0 + D_1^+(i)f_1 + D_2^+(i)f_2 + \dotsb.
\end{equation}
Similarly, the {\em preservation operator} $a^0(i)$ and
{\em annihilation operator} $a^-(i)$ are defined by
\begin{align}
a^0(i)f &= D_0^0(i)f_0 + D_1^0(i)f_1 + D_2^0(i)f_2 + \dotsb, \\
a^-(i)f &= D_0^-(i)f_0 + D_1^-(i)f_1 + D_2^-(i)f_2 + \dotsb,
\end{align}
respectively.
The operators $a^+(i), a^0(i), a^-(i)$, $i\in\{1,2,\dots,d\}$,
are generally called the {\em joint quantum operators}
associated to the random variables $X_1,X_2, \dots,X_d$.
Now Lemma \ref{lemma_restricted_quantum_decomposition} can be 
written as follows.

\begin{lemma}\label{lemma_quantum_decomposition}
For all $i \in \{1,2,\dots,d\}$ we have
\begin{equation}
X_i = a^-(i) + a^0(i) + a^+(i),
\end{equation}
where the domain of these operators is taken to be 
the space $F$ of all polynomial random variables
$f(X_1, X_2, \dots, X_d)$.
\end{lemma}

We have thus defined the quantum operators, 
first separately on each homogenous chaos space
$G_n$ for $n \geq 0$, 
and then extended them in a linear way
to the space $F$ of all polynomial random variables of
$X_1,X_2,\dots,X_d$.
We now mention an alternative way to define
the quantum operators directly on $F$.

For each $n \geq 0$ let $\bar{P}_n$ denote
the projection operator of $L^2(\Omega$, $P)$ onto $G_n$,
and $P_n$ the restriction of $\bar{P}_n$ to $F$.
The identity operator of $F$ is denoted by by $I$.
It follows from Lemma \ref{lemma_orthogonal} that
$P_mX_iP_n = 0$ for any non-negative integers $m$ and $n$
such that $|m-n|\geq 2$ and $1 \leq i \leq d$.
Then, for each $1 \leq i \leq d$ we have
\begin{align}
X_i &= IX_iI 
\nonumber\\
&= \bigg(\sum_{m=0}^{\infty}P_m\bigg)
   X_i\bigg(\sum_{n = 0}^{\infty}P_n\bigg) 
\nonumber\\
&= \sum_{|m - n| \leq 1}P_mX_iP_n
\nonumber\\
&= \sum_{n=0}^{\infty}P_{n+1}X_iP_n 
 + \sum_{n=0}^{\infty}P_nX_iP_n 
 + \sum_{n=1}^{\infty}P_nX_iP_{n-1}.
\end{align}
Each sum of the last expression corresponds to 
a quantum operator of index $i$, namely,
\begin{align}
a^+(i) &= \sum_{n = 0}^{\infty}P_{n+1}X_iP_n,
\label{02eqn:a+} \\
a^0(i) &= \sum_{n = 0}^{\infty}P_nX_iP_n,
\label{02eqn:a0} \\
a^-(i) &= \sum_{n = 1}^{\infty}P_nX_iP_{n-1}.
\label{02eqn:a-}
\end{align}
Moreover, from \eqref{02eqn:a+} and \eqref{02eqn:a-} we see that
\begin{equation}
(a^+(i))^{\star}
=\sum_{n=0}^{\infty}\left(P_{n+1}X_iP_n\right)^{\star}
=\sum_{n=0}^{\infty}P_nX_iP_{n+1}
=a^-(i)
\end{equation}
and from \eqref{02eqn:a0}
\begin{equation}
(a^0(i))^{\star}
=\sum_{n=0}^{\infty}\left(P_nX_iP_n\right)^{\star}
=\sum_{n=0}^{\infty}P_nX_iP_n
=a^0(i).
\end{equation}
The above formulas must be interpreted carefully,
since we speak about adjoints of unbounded operators.
To be cautious, 
we say that $a^-(i)$ is the polynomial dual of $a^+(i)$, 
while $a^0(i)$ is polynomially self-dual, 
which means that for any pair of polynomial random variables
$f,g \in F$ we have
\begin{equation}
\langle a^+(i)f, g\rangle = \langle f, a^-(i)g\rangle
\end{equation}
and
\begin{equation}
\langle a^0(i)f, g\rangle = \langle f, a^0(i)g\rangle,
\end{equation}
where $\langle \cdot,\cdot \rangle$ denotes
the inner product in $L^2(\Omega, P)$.

For all $i \in \{1,2,\dots,d\}$ 
we define the {\em semi-annihilation operator} of index $i$ by
\begin{equation}
U_i = a^-(i) + \frac{1}{2}\,a^0(i)
\end{equation}
and the {\em semi-creation operator} of index $i$ by
\begin{equation}
V_i = a^+(i) + \frac{1}{2}\,a^0(i).
\end{equation}
The following assertion is immediate from
Lemma \ref{lemma_quantum_decomposition}.

\begin{lemma}
For all $i \in \{1,2,\dots,d\}$ we have
\begin{equation}
X_i = U_i + V_i.
\end{equation}
\end{lemma}

Since $a^+(i)$ is the polynomial dual of $a^-(i)$
and $a^0(i)$ is polynomial self-dual,
we conclude that $V_i$ is the polynomial dual of $U_i$
for all $i \in \{1,2,\dots,d\}$.
We call $U_i$ and $V_i$, $i\in \{1,2,\dots,d\}$,
the {\em joint semi-quantum operators}
associated to $X_1,X_2,\dots,X_d$.

Another very important operator is the {\em number operator}
$N : F \to F$ defined by
\begin{equation}
N = \sum_{n = 0}^{\infty}nP_n,
\end{equation}
or equivalently by
\begin{equation}
Nf = nf,
\qquad f \in G_n, \quad n\in {\mathbb N} \cup \{0\}.
\end{equation}
%
%Finally, for every complex constant $c$
%we define the {\em second quantization operator of} $cI$, 
%where $I$ %denotes
%the identity operator of $F$, as:
%\begin{eqnarray}
%\Gamma(cI) & := & c^N \\
%& = & \sum_{n = 0}^{\infty}c^nP_n.
%\end{eqnarray}
%That means, for all $n \geq 0$, and all $f \in G_n$, we define:
%\begin{eqnarray}
%\Gamma(cI)f & := & c^nf.
%\end{eqnarray}

When $d=1$, we start with just random variable $X=X_1$
and, omitting the subscript $1$ we write
$a^+$, $a^0$, $a^-$, $V$ and $U$
for the creation, preservation, annihilation, semi-creation and
semi-annihilation operators, respectively.
In that case $F_n=F_{n-1}+\mathbb{C}x^n$
for all $n \geq 0$ and the codimension of $F_{n-1}$
in $F_n$ is at most $1$. 
Thus, the dimension of $G_n=F_n \ominus F_{n-1}$ is at most $1$
for all $n \geq 0$.
Consequently, if $G_n \neq 0$, 
then there exists a unique polynomial random variable $f_n(X)$
in $G_n$, having the leading coefficient (the coefficient of $X^n$)
equal to $1$. 
Due to this fact, the quantum decomposition $X=a^++a^0+a^-$ becomes
the classical recursive relation among monic orthogonal polynomials,
namely, for all $n \geq 0$ there exists real numbers $\alpha_n$
and $\omega_n$ such that
\begin{equation}
Xf_n(X)= f_{n+1}(X) + \alpha_nf_n(X) + \omega_nf_{n-1}(X).
\end{equation}
The numbers $\{\alpha_n\}_{n \geq 0}$ and $\{\omega_n\}_{n \geq 1}$
($\omega_0$ is not defined since $f_{-1}(X)= 0$)
are called the {\em Szeg\H{o}-Jacobi parameters} of $X$.
We thus come to
\begin{align}
a^+f_n(X) &= f_{n + 1}(X),
\\
a^0f_n(X) &= \alpha_nf_n(X),
\\
a^-f_n(X) &= \omega_nf_{n-1}(X).
\end{align}

\begin{remark}\normalfont
A similar idea of semi-quantum operators in case of $d=1$
appeared in \cite{hot2001}, where the adjacency matrix
of a Hamming graph was decomposed into a sum of 
``semi-quantum operators'' and its asymptotic spectral 
distribution was determined.
Later the ``semi-quantum operators'' was replaced with
the ``full'' quantum decomposition, see e.g., \cite{ho2007}.
\end{remark}

\section{Universal Commutation Relationships}

In what follows, for any pair of linear operators $A$ and $B$
we define the commutator by
\begin{equation}
[A, B] = AB - BA,
\end{equation}
provided both compositions $AB$ and $BA$ make sense.

It is very important to understand the commutation relations
satisfied by the quantum, semi-quantum and number operators
associated to a given family of random variables 
$X_1, X_2, \dots, X_d$ having finite moments of all orders.
Some of those commutation relationships are specific 
to a particular family random variables $X_1, X_2, \dots, X_d$,
and the first moments determine completely the random variables
in the sense of joint moments.
On the other hand, there is a set of commutation relationships 
that is common for all finite families of random variables.
These are called the {\em universal commutation relationships}
\cite{dps2020}. 

\begin{theorem}[Universal Commutators Theorem \cite{dps2020}]
\label{universal_commutators_theorem} 
Let $X_1, X_2, \dots, X_d$ be random variables defined on 
a common probability space and having finite moments of all orders.
Let $a^+(i), a^0(i), a^-(i)$, $V_i, U_i$ and $N$ be their 
joint creation, preservation, annihilation, semi-creation,
semi-annihilation, and number operators, where $i\in\{1,2,\dots,d\}$.
Then, for all $i \in \{1,2,\dots,d\}$ we have
\begin{align}
&[N, a^+(i)] = a^+(i), 
\label{Na+} \\
&[N, a^0(i)] = 0, 
\label{Na0} \\
&[a^-(i), N] = a^-(i), 
\label{a-N} \\
&[N, V_i] = a^+(i), 
\label{NV_i} \\
&[U_i, N] = a^-(i), 
\label{U_iN} \\
&[N, X_i] =V_i - U_i=a^+(i) - a^-(i). 
\label{NX_i}
\end{align}
\end{theorem}

The above commutation relationships can be combined
with the rules of ``commutative probability", 
i.e., $X_iX_j = X_jX_i$ for all $i,j \in \{1,2,\dots,d\}$.
These rules, established in \cite{aks04,aks07,ps14},
are as follows: for all $i,j\in\{1,2,\dots,d\}$ we have
\begin{align}
&[a^+(i), a^+(j)] = 0,
\\
&[a^0(i), a^+(j)] = [a^0(j), a^+(i)],
\\
&[a^-(i),a^+(j)] + [a^0(i),a^0(j)] + [a^+(i),a^-(j)] = 0,
\\
&[a^-(i),a^0(j)] = [a^-(j), a^0(i)],
\\
&[a^-(i), a^-(j)] = 0,
\\
&[X_i, V_j] = [X_j, V_i],
\\
&[U_i, X_j] = [U_j, X_i],
\\
&[X_i, V_j]^{\star} = [X_i, V_j],
\\
&[U_i, X_j]^{\star} = [U_i, X_j].
\end{align}

\section{Position-Momentum Decomposition of Linear Operators 
Defined on $F$}

Let $\mathbb{C}[x_1,x_2,\dots,x_d]$ be the algebra of polynomials
with complex coefficients in $d$ variables $x_1, x_2, \dots, x_d$.

\begin{definition}\normalfont
For any $\vec{\alpha}=(\alpha_1,\alpha_2,\dots,\alpha_d)
\in (\mathbb{N} \cup \{0\})^d$, we define the
{\em $\vec{\alpha}$-position operator} 
$X^{\vec{\alpha}}:\mathbb{C}[x_1,x_2, \dots,x_d]\to
\mathbb{C}[x_1,x_2,\dots,x_d]$ by
\begin{equation}
X^{\vec{\alpha}}f\left(x_1, x_2, \dots, x_d\right) 
=x_1^{\alpha_1}x_2^{\alpha_2} \cdots x_d^{\alpha_d}
 f(x_1, x_2, \dots, x_d),
\end{equation}
and {\em $\vec{\alpha}$-momentum operator} 
$D^{\vec{\alpha}} : \mathbb{C}[x_1, x_2, \dots, x_d]
\to \mathbb{C}[x_1, x_2, \dots, x_d]$ by
\begin{equation}
D^{\vec{\alpha}}f\left(x_1, x_2, \dots, x_d\right)
= \frac{\partial^{\alpha_1}}{\partial x_1^{\alpha_1}}\frac{\partial^{\alpha_2}}{\partial x_2^{\alpha_2}}
\cdots \frac{\partial^{\alpha_d}}{\partial x_d^{\alpha_d}}
f(x_1, x_2, \dots, x_d).
\end{equation}
\end{definition}

The position operator $X^{\vec{\alpha}}$ becomes naturally
a multiplication operator on the space $F$ of polynomial random 
variables as a composition of the multiplication
operators by $X_1, X_2, \dots, X_d$.
On the other hand, the momentum operator $D^{\vec{\alpha}}$
may not be well-defined as an operator from $F$ into itself.
In fact, if the joint probability distribution $\mu$ of 
$X_1, X_2, \dots, X_d$ is supported on the algebraic variety 
$V = \{\vec{x} \in {\mathbb R}^d \mid f_0(\vec{x}) = 0\}$ 
for some polynomial $f_0(x_1$, $x_2$, $\dots$, $x_d) 
\in \mathbb{C}[x_1, x_2, \dots,x_d]$, 
then we have $g(X_1,X_2,\dots,X_d)=h(X_1,X_2,\dots,X_d)$ 
holds almost surely for any pair of polynomials $g$ and $h$
such that $f_0 \mid (g - h)$, 
but $D^{\vec{\alpha}}g(X_1, X_2, \dots, X_d)$ may not 
be almost surely equal to 
$D^{\vec{\alpha}}h(X_1,X_2,\dots, X_d)$
since $f_0 \nmid (D^{\vec{\alpha}}g - D^{\vec{\alpha}}h)$ may happen.

\begin{definition}\normalfont
Let $k$ be a fixed integer and $T:F \to F$ a linear operator.
We say that $T$ is {\em $k$-faithful to the gradation
$\{F_n\}_{n \geq 0}$} if
\begin{equation}
TF_n \subseteq F_{n + k},
\qquad n\in\mathbb{N}\cup\{0\}.
\end{equation}
Here we tacitly understand that $F_{n+k}=\{0\}$ if $n+k<0$.
\end{definition}

\begin{theorem}[\cite{dps2020}]
\label{X_D_decomposition}
If $T$ is $k$-faithful, 
then there exists a unique family of polynomial random variables
$\{A_{\vec{\alpha}}(X_1,X_2,\dots,X_d)\}_{\vec{\alpha}
\in (\mathbb{N}\cup\{0\})^d}$ with
\begin{equation}
\deg(A_{\vec{\alpha}}) \leq |\vec{\alpha}|+ k,
\qquad
\vec{\alpha}
=(\alpha_1, \alpha_2, \dots, \alpha_d)\in ({\mathbb N} \cup \{0\})^d,
\end{equation}
where $|\vec{\alpha}|= \alpha_1 + \alpha_2 + \dotsb + \alpha_d$ is
the length of $\vec{\alpha}$ and $\deg(A_{\vec{\alpha}})$
denotes the degree of $A_{\vec{\alpha}}$,
such that 
\begin{equation}\label{explicit_X_D_decomposition}
T(f(X_1, X_2, \dots, X_d))
=\sum_{\vec{\alpha} \in ({\mathbb N} \cup \{0\})^d}
A_{\vec{\alpha}}(X_1, X_2, \dots, X_d)
(D^{\vec{\alpha}}f)(X_1, X_2, \dots, X_d)
\end{equation}
for all $f(x_1, x_2,\dots, x_d)\in\mathbb{C}[x_1,x_2,\dots,x_d]$.
The uniqueness of $\{A_{\vec{\alpha}}
(X_1,X_2,\dots,X_d\}_{\alpha\in(\mathbb{N}\cup\{0\})^d}$ is
understood in the sence of $\mu$-almost sure, 
where $\mu$ is the joint probability distribution of
$X_1, X_2, \dots, X_d$.
\end{theorem}

For convenience formula \eqref{explicit_X_D_decomposition} is
written as
\begin{equation}
T =\sum_{\vec{\alpha} \in {{\mathbb N} \cup \{0\}}^d}
A_{\vec{\alpha}}D^{\vec{\alpha}},
\end{equation}
which is called the {\em position-momentum decomposition} of $T$.

\begin{remark}\normalfont
For a given polynomial random variable $Y \in F$
there is ambiguity of choosing a polynomial $f$ 
such that $Y = f(X_1, X_2, \dots, X_d)$, $\mu$-a.s.
Once such a polynomial $f$ is chosen the same polynomial
must be used for all the terms in the right-hand side of 
formula \eqref{explicit_X_D_decomposition}).
\end{remark}

For $\vec{\alpha}=(\alpha_1,\dots,\alpha_d),
\vec{\beta}=(\beta_1,\dots,\beta_d)\in(\mathbb{N}\cup\{0\})^d$
we write $\vec{\alpha} \preceq \vec{\beta}$ if 
$\alpha_i \leq \beta_i$ for all $1 \leq i \leq d$.
We then write $\vec{\alpha} \prec \vec{\beta}$ 
if $\vec{\alpha} \preceq \vec{\beta}$
and $\vec{\alpha} \neq \vec{\beta}$.
For a polynomial $f(x_1, x_2, \dots, x_d)$ let
$m_i$ be the highest power at which $x_i$ appears in $f$
and set $\vec{m}= (m_1,m_2,\dots, m_d) \in (\mathbb{N}\cup\{0\})^d$.
Then formula \eqref{explicit_X_D_decomposition} becomes
\begin{equation}\label{04eqn:explicit_X_D_decomposition}
T(f\left(X_1, X_2, \dots, X_d\right))
=\sum_{\vec{\alpha} \preceq \vec{m}}
A_{\vec{\alpha}}(X_1, X_2, \dots, X_d)
(D^{\vec{\alpha}}f)(X_1, X_2, \dots, X_d).
\end{equation}
In fact, if we differentiate $f$ with respect to
one of the variables $x_1, x_2, \dots, x_d$ more times than
its highest power, we get $0$.
Thus \eqref{explicit_X_D_decomposition}
is reduced to \eqref{04eqn:explicit_X_D_decomposition}.
The above argument is based on the form of
\eqref{04eqn:explicit_X_D_decomposition},
where the ``position factors" contained 
in the multiplication operators
$A_{\vec{\alpha}}(\vec{X}^{\vec{\alpha}})$ 
is placed to the left of all the ``momentum factors" 
$D^{\vec{\alpha}}$.
We also note that
a mixed composition of position and momentum operators
$x_1$, $\dots$, $x_d$, $\partial/\partial x_1$, \dots, $\partial x_d$, that do not obey the left-position right-momentum order
is always rearranged by using the famous Heisenberg 
commutation relations:
\begin{equation}
\left[\frac{\partial}{\partial x_i}, x_j\right] = \delta_{i,j}I,
\qquad
1 \leq i,j \leq d,
\end{equation}
where $\delta_{i,j}$ denotes the Kronecker symbol
and I the identity operator. 

For more detailed information on
the position-momentum decomposition,
see \cite{dps2020}.

\section{Position-Momentum Decomposition of Meixner Random Variables}

We start with the classical Meixner random variables.

\begin{definition}\normalfont
A real-valued random variable $X$ having finite moments of all orders
is called a {\em classical Meixner random variable} if
its Szeg\H{o}-Jacobi parameters are of the form:
\begin{equation}
\alpha_n = \alpha n + \alpha_0,
\qquad
\omega_n = \beta n^2 + (t - \beta)n,
\qquad
n \geq 1,
\end{equation}
where $\alpha$, $\alpha_0$, $\beta$, and $t$ are real constants
satisfying one of the following conditions:
\begin{enumerate}
\item $\beta \geq 0$ and $t \geq 0$;
\item $\beta<0$ and $t \in -{\mathbb N}\beta$.
\end{enumerate}
\end{definition}

If $t = 0$, then $\omega_1 = 0$ and hence $X$ becomes
a constant random variable.
This case being not interesting, we only consider the case of $t>0$.
We may also assume that $\alpha \geq 0$ 
since otherwise, we can replace $X$ by $-X$.

We will follow the following steps:

\noindent
{\bfseries Step 1.} Find the commutator $[U$, $X]$ between 
the semi-annihilation operator $U$ and multiplication operator by $X$. 
For all $n \in {\mathbb N}$, 
if $f_n(X)$ denotes the orthogonal polynomial of degree $n$
with leading coefficient equal to $1$, we have
\begin{align}
[U, X]f_n(X)
&= [U, U + V]f_n(X) \nonumber\\
&= [U, V]f_n(X) \nonumber\\
&= \left[a^- + \frac{1}{2}a^0, a^+ + \frac{1}{2}a^0\right]f_n(X) 
\nonumber\\
&= \left[a^-, a^+\right]f_n(X) 
 + \frac{1}{2}\left[a^-, a^0\right]f_n(X) 
 + \frac{1}{2}\left[a^0, a^+\right]f_n(X) 
\nonumber\\
&= a^-a^+f_n(X) - a^+a^-f_n(X) 
 + \frac{1}{2}\left[a^-a^0f_n(X) - a^0a^-f_n(X)\right] 
\nonumber\\
&\qquad + \frac{1}{2}\left[a^0a^+f_n(X) - a^+a^0f_n(X)\right]
\nonumber\\
&= a^-f_{n + 1}(X) - a^+\omega_nf_{n-1}(X) 
 + \frac{1}{2}\left[a^-\alpha_nf_n(X) - a^0\omega_nf_{n-1}(X)\right] 
\nonumber\\
&\qquad +\frac{1}{2}\left[a^0f_{n+1}(X) - a^+\alpha_nf_{n-1}(X)\right] 
\nonumber\\
&= \left(\omega_{n+1} - \omega_n\right)f_n(X) 
 + \frac{1}{2}\omega_n\left(\alpha_n - \alpha_{n-1}\right)f_{n-1}(X)
 + \frac{1}{2}\left(\alpha_{n+1} - \alpha_n\right)f_{n+1}(X). 
\nonumber
\end{align}
For $n \geq 1$ we have
\begin{align}
\alpha_n - \alpha_{n-1} 
&=\alpha(n+1) + \alpha_0 - \alpha n - \alpha_0
=\alpha, 
\\
\omega_{n+1} - \omega_n 
&= \beta(n+1)^2 + (t-\beta)(n+1) - \beta n^2 - (t-\beta)n
= 2\beta n + t.
\end{align}
Thus, for all $n \geq 1$ we obtain:
\begin{align}
[U, X]f_n(X) 
&= (2\beta n+t)f_n(X) + \frac{1}{2}\alpha\omega_nf_{n-1}(X) 
 + \frac{1}{2}\alpha f_{n+1}(X) 
\nonumber\\
&= (2\beta N+tI)f_n(X) + \frac{\alpha}{2}a^-f_n(X) 
 + \frac{\alpha}{2}a^+f_n(X) 
\nonumber\\
&= (2\beta N+tI)f_n(X) + \frac{\alpha}{2}(a^- + a^+)f_n(X)
\nonumber\\
&= (2\beta N + tI)f_n(X) + \frac{\alpha}{2}(X-a^0)f_n(X)
\nonumber\\
&= (2\beta N+tI)f_n(X) 
+ \frac{\alpha}{2}(X - \alpha N - \alpha_0I)f_n(X).
\nonumber
\end{align}
The last formula is also true for $f_0(X) = 1$ since
\begin{align}
[U, X]f_0(X) 
&= UXf_0(X) - XUf_0(X) 
\nonumber\\
&= \left(\frac{1}{2}a^0 + a^-\right)
 \left[f_1(X) + \alpha_0f_0(X)\right] 
 - X\left(\frac{1}{2}a^0 + a^-\right)f_0(X) 
\nonumber\\
&= \frac{\alpha_1}{2}f_1(X) + \frac{\alpha_0^2}{2}f_0(X) 
 + \omega_1f_0(X) - X\frac{\alpha_0}{2}f_0(X) 
\nonumber\\
&= \frac{\alpha_1}{2}f_1(X) + \frac{\alpha_0^2}{2}f_0(X) 
 + \omega_1f_0(X) - \frac{\alpha_0}{2}(f_1(X) + \alpha_0f_0(X))
\nonumber\\
&= \frac{\alpha}{2}f_1(X) + \omega_1f_0(X) 
\nonumber\\
&= tf_0(X) + \frac{\alpha}{2}\left(X - \alpha_0\right)
\nonumber
\end{align}
and $Nf_0(X) = 0f_0(X) = 0$.
Thus, we obtain
\begin{equation}\label{04eqn:[U,X]}
[U, X] 
= \frac{\alpha}{2}X + \frac{4\beta - \alpha^2}{2}N 
+ \frac{2t - \alpha\alpha_0}{2}I.
\end{equation}
We set
\begin{equation}
\Delta= \alpha^2 - 4\beta, 
\qquad
\tau = 2t - \alpha\alpha_0.
\end{equation}
Then \eqref{04eqn:[U,X]} becomes
\begin{equation}
[U, X] 
= \frac{\alpha}{2}X - \frac{\Delta}{2}N + \frac{\tau}{2}I.
\end{equation}
We consider two cases.

\noindent
{\bfseries Case 1.} $\Delta = 0$.
We obtain
\begin{equation}\label{comm_U_X_1_Meixner}
[U, X] = \frac{\alpha}{2}X + \frac{\tau}{2}I, 
\end{equation}
which means that the commutator $[U, X]$ is a linear combination of 
$X$ and $I$. 
After \cite{ps14} we call such a random variable is called 
{\em $1$-Meixner random variable}.

\noindent
{\bfseries Case 2.} $\Delta \neq 0$.
Using the universal commutator theorem, 
we compute the following double commutator: 
\begin{align}
[[U, X], X] 
&= \left[\frac{\alpha}{2}X - \frac{\Delta}{2}N 
 + \frac{\tau}{2}I, X\right] 
\nonumber\\
&= -\frac{\Delta}{2}[N, X] 
\nonumber\\
&= -\frac{\Delta}{2}\left(V - U\right) 
\nonumber\\
&= -\frac{\Delta}{2}\left(X - 2U\right). 
\label{comm_U_X_2_Meixner}
\end{align}
After \cite{ps14} such a random variable is called
{\em $2$-Meixner random variables}.

\medskip
\noindent
{\bfseries Step 2.} 
Find the position-momentum decomposition of
the semi-quantum operators $U$ and $V$.
Indeed, let the position-momentum decomposition of $U$ be
\begin{equation}
U = \sum_{n = 0}^{\infty}A_n(X)D^n,
\end{equation}
where for all $n \in {\mathbb N} \cup \{0\}$, 
$A_n$ is a polynomial of degree at most $n$, 
since $U$ maps $F_n$ into $F_n$.

\medskip
\noindent
{\bf Case 1.} $\Delta = 0$.
Substituting $U$ by its position-momentum decomposition in 
\eqref{comm_U_X_1_Meixner}, we obtain:
\begin{equation}
\frac{\alpha}{2}X + \frac{\tau}{2}I 
= [U, X]
= \sum_{n = 0}^{\infty}\left[A_n(X)D^n, X\right].
\end{equation}
Here we can interchange the sum with the commutator, 
since we apply the operators $A_n(X)D^n$ only
to polynomial random variables $f(X) \in F$, 
and for each fixed $f$, for all $n > N$, 
where $N$ is the degree of $f$, we have $D^n = 0$. 
Thus after a while we are adding only zero terms.
Applying the Leibniz commutator rule:
\begin{equation}
[B_1B_2 \cdots B_k, C] 
=\sum_{i=1}^kB_1 \cdots B_{i-1}[B_i, C]B_{i+1} \cdots B_k,
\end{equation}
we obtain:
\begin{align}
\frac{\alpha}{2}X + \frac{\tau}{2}I 
&= \sum_{n = 0}^{\infty}\left[A_n(X)D^n, X\right] \nonumber\\
&= \sum_{n = 1}^{\infty}
 \left\{[A_n(X), X]D^n 
 +\sum_{i=1}^n A_n(X) D^{i-1}[D,X]D^{n-i}\right\}
\nonumber\\
&= \sum_{n=1}^{\infty}
 \left\{0D^n + \sum_{i=1}^nA_n(X)D^{i-1}ID^{n-i}\right\}
\nonumber\\
&= \sum_{n=1}^{\infty}nA_n(X)D^{n-1}
\nonumber\\
&= \sum_{n=0}^{\infty}(n+1)A_{n+1}(X)D^n.
\end{align}
We can write the above formula as
\begin{equation}
\left(\frac{\alpha}{2}X 
 + \frac{\tau}{2}I\right)D^0 + 0D^1 + 0D^2 + \cdots 
= A_1(X)D^0 + 2A_2(X)D^1 + 3A_3(X)D^3 + \cdots.
\end{equation}
Using the uniqueness of the position-momentum decomposition,
we can equate the left-position coefficients of 
the corresponding right-momentum factors $D^n$ 
from the two sides of the above equality, and obtain:
\begin{align}
A_1(X) &= \frac{\alpha}{2}X + \frac{\tau}{2}I, \\
nA_n(X) &= 0, \qquad n \geq 2.
\end{align}
The left-position coefficient $A_0(X)$ can be 
obtained from the fact that, 
for the constant polynomial $f_0 = 1$ we have:
\begin{equation}
U1=\left(a^- + \frac{1}{2}a^0\right)1 \nonumber\\
=\frac{\alpha_0}{2}1
\end{equation}
and
\begin{equation}
U1 = \sum_{n = 0}^\infty A_n(X)D^n1
=A_0(X)D^01
=A_0(X)1.
\end{equation}
Thus, we obtain:
\begin{equation}
A_0(X)=\frac{\alpha_0}{2}.
\end{equation}
Therefore, the position-momentum quantum decomposition of the 
$1$-Meixner random variable with parameters $\alpha$, $\alpha_0$, 
and $t$ (remember that $\beta =\alpha^2/4$), is
\begin{equation}\label{U_1_Meixner_decomp}
U = \frac{\alpha_0}{2}I + \frac{1}{2}(\alpha X + \tau)D. 
\end{equation}
In view of $V = X - U$, 
we conclude that the position-momentum decomposition of $V$ is
given by
\begin{equation}\label{V_1_Meixner}
V=\left(X - \frac{\alpha_0}{2}\right)I 
 -\frac{1}{2}(\alpha X + \tau I)D. 
\end{equation}

\medskip
\noindent
{\bfseries Case 2.} $\Delta \neq 0$.
Substituting $U$ by its position-momentum decomposition
in \eqref{comm_U_X_2_Meixner}, we obtain
\begin{equation}
\left[\left[\sum_{n = 0}^{\infty}A_n(X)D^n, X\right], X\right] 
= -\frac{\Delta}{2}\left(X - 2\sum_{n = 0}^{\infty}A_n(X)D^n\right).
\end{equation}
Using the Leibniz commutator formula, we obtain
\begin{equation}
\sum_{n = 1}^{\infty}\left[nA_n(X)D^{n - 1}, X\right] 
=-\frac{\Delta}{2}\left\{X - 2A_0(X)\right\}I 
 +\Delta\sum_{n = 1}^{\infty}A_n(X)D^n.
\end{equation}
Using the Leibniz commutator formula again, we obtain
\begin{equation}
\sum_{n = 2}^{\infty}n(n - 1)A_n(X)D^{n - 2} 
= -\frac{\Delta}{2}\left\{X - 2A_0(X)\right\}I 
 +\Delta\sum_{n = 1}^{\infty}A_n(X)D^n,
\end{equation}
which, after making the change of variable $n \mapsto n-2$
in the left sum, becomes:
\begin{equation}
\sum_{n=0}^{\infty}(n+2)(n+1)A_{n+2}(X)D^n
=-\frac{\Delta}{2}\left\{X - 2A_0(X)\right\}I 
 +\Delta\sum_{n=1}^{\infty}A_n(X)D^n.
\end{equation}
Equating the left-position coefficients of the corresponding
right-momentum factors $D^n$, for $n\geq 0$
we obtain the recursive formula:
\begin{equation}\label{05eqn:recursive formula for A}
A_{n+2}(X) 
=\frac{\Delta}{(n+2)(n+1)}
 \left(A_n(X)-\frac{1}{2}\delta_{n,0}X\right),
\end{equation}
where $\delta_{p,q}$ denotes the Kronecker symbol.
Iterating the recursive relation \eqref{05eqn:recursive formula for A},
we obtain
\begin{align}
A_{2n}(X) 
&= \frac{\Delta^n}{(2n)!}\left(A_0(X) - \frac{1}{2}X\right) \label{even_recursive} \\
A_{2n+1}(X) 
&= \frac{\Delta^n}{(2n + 1)!}A_1(X),
\label{odd_recursive}
\end{align}
for all $n \in \mathbb{N}$.
The left-position coefficients $A_0(X)$ and $A_1(X)$ can be
found by applying $U$ to $f_0(X) = 1$ and
$f_1(X) = X - \alpha_0$. 
We have, as before
\begin{equation}
U1
= \frac{1}{2}a^01
= \frac{\alpha_0}{2}1
\end{equation}
and
\begin{equation}
U1 = A_0(X)1,
\end{equation}
which means that $A_0(X) = \alpha_0/2$.
Moreover, we have
\begin{align}
Uf_1(X) 
&= a^-f_1(X) + \frac{1}{2}a^0f_1(X) \nonumber\\
&= \omega_1f_0(X) + \frac{\alpha_1}{2}f_1(X) \nonumber\\
&= (\beta \cdot 1^2 +(t-\beta)1)1 + \frac{\alpha\cdot 1+\alpha_0}{2}
 (X - \alpha_0) \nonumber\\
&= \frac{\alpha + \alpha_0}{2}(X-\alpha_0)+t
\label{05eqn:Uf_1(X)1}
\end{align}
and also
\begin{align}
Uf_1(X) 
&= (A_0(X)I + A_1(X)D)(X - \alpha_0) 
\nonumber \\
&= \frac{\alpha_0}{2}(X - \alpha_0) + A_1(X).
\label{05eqn:Uf_1(X)2}
\end{align}
From \eqref{05eqn:Uf_1(X)1} and \eqref{05eqn:Uf_1(X)2} we see that
\begin{equation}
A_1(X) = \frac{\alpha}{2}\left(X - \alpha_0\right) + t.
\end{equation}
It follows from \eqref{even_recursive} and \eqref{odd_recursive} 
that for all $n \in {\mathbb N}$ we have
\begin{align}
A_{2n}(X) 
&= -\frac{1}{2}\cdot\frac{\Delta^n}{(2n)!}(X - \alpha_0)
\\
A_{2n+1}(X) 
&= \frac{\Delta^n}{(2n+1)!}
 \left[\frac{\alpha}{2}(X - \alpha_0)+t\right].
\end{align}
Let $\delta$ be a complex square root of $\Delta$, 
that means, $\Delta = \delta^2$. 
Since $\Delta \neq 0$, we have $\delta \neq 0$.
Thus we obtain the position-momentum decomposition of $U$ as
\begin{align}
U &=\sum_{n=0}^{\infty}A_{2n}(X)D^{2n} 
 +\sum_{n=0}^{\infty}A_{2n+1}(X)D^{2n+1}
\nonumber\\
&=\frac{\alpha_0}{2}I 
 -\frac{1}{2}(X-\alpha_0)
  \sum_{n=1}^{\infty}\frac{\Delta^n}{(2n)!}D^{2n} 
 +\left[\frac{\alpha}{2}(X-\alpha_0)+t\right]
  \sum_{n=0}^{\infty}\frac{\Delta^n}{(2n+1)!}D^{2n+1}.
\label{U_formula_before_translation}
\end{align}

With each complex number $c$ we associate the translation operator
$\mathcal{T}_c: F \to F$ defined by
\begin{equation}
\mathcal{T}_cf(X) = f(X + c).
\end{equation}
We then see by the Taylor formula that the position-momentum 
decomposition of $\mathcal{T}_{\delta}$ is given by
\begin{equation}\label{Taylor_delta}
\mathcal{T}_{\delta}
=\sum_{n=0}^{\infty}\frac{\delta^n}{n!}D^n 
=\exp(\delta D) 
\end{equation}
and similarly,
\begin{equation}\label{Taylor_minus_delta}
\mathcal{T}_{-\delta} 
= \sum_{n=0}^{\infty}\frac{(-\delta)^n}{n!}D^n
= \exp(-\delta D). 
\end{equation}
Adding and subtracting both sides of \eqref{Taylor_delta} and 
\eqref{Taylor_minus_delta}, we obtain
\begin{align}
\sum_{n=0}^{\infty}\frac{\delta^{2n}}{(2n)!}D^{2n} 
&=\frac{1}{2}(\mathcal{T}_\delta + \mathcal{T}_{-\delta})
= \cosh(\delta D), 
\label{cosh(D)} \\
\sum_{n=0}^{\infty}\frac{\delta^{2n+1}}{(2n+1)!}D^{2n+1} 
&=\frac{1}{2}(\mathcal{T}_\delta - \mathcal{T}_{-\delta})
= \sinh(\delta D). 
\label{sinh(D)}
\end{align}
In view of $\tau = 2t-\alpha\alpha_0$ 
formula \eqref{U_formula_before_translation} becomes 
\begin{equation}\label{U_2_Meixner_formula}
U = \frac{\alpha_0}{2}I 
 -\frac{1}{4}(X - \alpha_0)
  (\mathcal{T}_\delta + \mathcal{T}_{-\delta} - 2I)
 +\frac{1}{4\delta}(\alpha X + \tau)
  (\mathcal{T}_\delta-\mathcal{T}_{-\delta}).
\end{equation}
Let us observe that as $\Delta \to 0$, 
using operatorial limits in a weak sense, we have
\begin{align}
\lim_{\Delta \to 0}U 
&= \frac{\alpha_0}{2}I - \frac{1}{4}\left(X - \alpha_0\right)
\left(\lim_{\delta \to 0}\mathcal{T}_{\delta} + \lim_{\delta \to 0}\mathcal{T}_{-\delta} - 2I\right)
+ \frac{1}{2}(\alpha X + \tau)\lim_{\delta \to 0}\frac{\mathcal{T}_{\delta} - \mathcal{T}_{-\delta}}{2\delta} 
\nonumber\\
&= \frac{\alpha_0}{2}I 
 - \frac{1}{4}(X - \alpha_0)\cdot 0 + \frac{1}{2}(\alpha X+\tau)D 
\nonumber\\
&= \frac{\alpha_0}{2}I + \frac{1}{2}(\alpha X + \tau)D,
\nonumber
\end{align}
which is exactly the position-momentum decomposition 
for the semi-quantum operator of $1$-Meixner random variables 
found in \eqref{U_1_Meixner_decomp}. 
This shows that the $1$-Meixner random variables can be understood
as a limiting case, as $\Delta \to 0$, 
of the $2$-Meixner random variables.
The semi-creation operator $V$ can be found very easily now as follows:
\begin{align}
V &= X - U 
\nonumber\\
&= X - \frac{\alpha_0}{2}I + \frac{1}{4}(X - \alpha_0)
 (\mathcal{T}_{\delta} + \mathcal{T}_{-\delta} - 2I)
 - \frac{1}{4\delta}(\alpha X + \tau)
 (\mathcal{T}_{\delta} - \mathcal{T}_{-\delta}). \label{V_2_Meixner_formula}
\end{align}

\medskip
\noindent
{\bfseries Step 3.} 
Find the position-momentum decomposition of the number operator $N$.
Since $NF_n \subseteq F_n$ for all $n \geq 0$, 
the position-momentum decomposition of $N$ is of the form:
\begin{equation}
N = \sum_{n = 0}^{\infty}A_n(X)D^n,
\end{equation}
where $A_n(X) \in F_n$ for all $n \geq 0$.
Using the universal commutator formula:
\begin{equation}
[N, X] = X - 2U,
\end{equation}
and replacing $U$ by its position-momentum decomposition
which is found in Step 2, we obtain
\begin{align*}
&\left[\sum_{n = 0}^{\infty}A_n(X)D^n, X\right] \\
&\qquad=
X - \alpha_0I 
+(X-\alpha_0)\sum_{n=1}^{\infty}\frac{\delta^{2n}}{(2n)!}D^{2n}
-\frac{1}{\delta}(\alpha X+\tau)
 \sum_{n=0}^\infty\frac{\delta^{2n+1}}{(2n+1)!}D^{2n+1},
\end{align*}
with the understanding that if $\delta = 0$, 
then a weak limit must be taken to make sense of the above formula.
Using the Leibniz commutator formula, 
the above formula becomes
\begin{equation}
\sum_{n = 1}^{\infty}nA_n(X)D^{n - 1} 
=(X - \alpha_0)\sum_{n=0}^{\infty}
  \frac{\delta^{2n}}{(2n)!}D^{2n} 
 -\frac{1}{\delta}(\alpha X + \tau)
  \sum_{n=0}^{\infty}\frac{\delta^{2n+1}}{(2n+1)!}D^{2n+1},
\end{equation}
which implies that
\begin{equation}
\sum_{n=0}^{\infty}(n+1)A_{n+1}(X)D^n
=(X-\alpha_0)\sum_{n=0}^{\infty}
 \frac{\delta^{2n}}{(2n)!}D^{2n}
 -\frac{1}{\delta}(\alpha X+\tau)
  \sum_{n=0}^{\infty}\frac{\delta^{2n+1}}{(2n+1)!}D^{2n+1}.
\end{equation}
Equating the left-position coefficients of 
the corresponding right-momentum terms $D^n$ for $n \geq 0$, 
we obtain 
\begin{align}
A_{2n + 1} 
&= (X - \alpha_0)\frac{\delta^{2n}}{(2n + 1)!},
\qquad n \geq 0,
\\
A_{2n} 
&= -\frac{1}{\delta}\cdot\frac{\delta^{2n-1}}{(2n)!}(\alpha X+\tau),
\qquad n \geq 1.
\end{align}
We see from $N1 = 0$ that $A_0(X) = 0$ 
($N1 = A_0(X)1 + A_1(X)D1 + \cdots$).
Thus, the position-momentum decomposition of $N$ is given by
\begin{equation}\label{N_position_momentum_decomp}
N = \frac{X - \alpha_0}{\delta}
 \sum_{n=0}^{\infty}\frac{\delta^{2n+1}}{(2n+1)!}D^{2n+1}
 -\frac{1}{\delta^2}(\alpha X + \tau)
 \sum_{n=1}^{\infty}\frac{\delta^{2n}}{(2n)!}D^{2n}.
\end{equation}
This can also be written as
\begin{equation}\label{N_translation_formula}
N = \frac{X-\alpha_0}{2\delta}
 (\mathcal{T}_{\delta} - \mathcal{T}_{-\delta})
 -\frac{1}{2\delta^2}(\alpha X + \tau)
 (\mathcal{T}_\delta + \mathcal{T}_{-\delta} - 2I).
\end{equation}
In particular, for $\Delta = \delta^2 = 0$ we have
\begin{align}
N &=
(X - \alpha_0)
 \lim_{\delta \to 0}\left\{\frac{1}{2\delta}
  (\mathcal{T}_{\delta} - \mathcal{T}_{-\delta})\right\}
- \frac{1}{2}(\alpha X + \tau)
 \lim_{\delta \to 0}\left\{\frac{1}{\delta^2}
  (\mathcal{T}_{\delta} + \mathcal{T}_{-\delta} - 2I)\right\} \nonumber\\
&=(X-\alpha_0)D - \frac{1}{2}(\alpha X + \tau) D^2. \label{N_position_momentum_delta_zer0}
\end{align}

\medskip
\noindent
{\bfseries Step 4.} 
Find the position-momentum decomposition of the quantum operators
$a^+, a^-$ and $a^0$.
We can follow the strategy outlined in \cite{dps2020} and
find first $a^+$ and $a^-$ using the universal commutation formulas:
\begin{equation}
[N, V] = a^+,
\qquad
[U, N] = a^-,
\end{equation}
and then $a^0 = X - a^+ - a^-$. 
On the other hand, since $\alpha_n = \alpha n + \alpha_0$
for $n \geq 0$, we have
\begin{equation}
a^0 = \alpha N + \alpha_0I.
\end{equation}
Using $N$, we find $a^0$ immediately as
\begin{equation}\label{preservation_operator_formula}
a^0=\frac{\alpha}{2\delta}
 (X - \alpha_0)(\mathcal{T}_\delta - \mathcal{T}_{-\delta})
 -\frac{\alpha}{2\delta^2}(\alpha X + \tau)
 (\mathcal{T}_\delta + \mathcal{T}_{-\delta} - 2I) + \alpha_0I. 
\end{equation}
In view of $U = a^- + (1/2)a^0$ and $V = a^+ + (1/2)a^0$
we can easily find $a^-$ and $a^+$.
Thus, using formulas \eqref{U_2_Meixner_formula},
\eqref{V_2_Meixner_formula} and \eqref{preservation_operator_formula},
we obtain
\begin{equation}
a^- = U - \frac{1}{2}a^0.
\end{equation}
We have:
\begin{align}
a^- 
&=\frac{\alpha_0}{2}I 
 -\frac{1}{4}(X - \alpha_0)
   (\mathcal{T}_\delta + \mathcal{T}_{-\delta} - 2I)
 +\frac{1}{4\delta}(\alpha X + \tau)
   (\mathcal{T}_{\delta} - \mathcal{T}_{-\delta})
\nonumber\\
&\qquad - \frac{\alpha}{4\delta}
   (X - \alpha_0)(\mathcal{T}_{\delta} - \mathcal{T}_{-\delta}) 
 +\frac{\alpha}{4\delta^2}(\alpha X + \tau)
   (\mathcal{T}_{\delta} + \mathcal{T}_{-\delta} - 2I) 
 -\frac{\alpha_0}{2}I
\nonumber\\
&= \frac{\tau + \alpha\alpha_0}{4\delta}
   (\mathcal{T}_{\delta} - \mathcal{T}_{-\delta})
 +\frac{\alpha_0\delta^2 + \alpha\tau}{4\delta^2}
   (\mathcal{T}_{\delta} + \mathcal{T}_{-\delta} - 2I)
\nonumber\\
&\qquad + \frac{\alpha^2 - \delta^2}{4\delta^2}X
   (\mathcal{T}_\delta + \mathcal{T}_{-\delta} - 2I).
\nonumber
\end{align}
Since $\tau + \alpha\alpha_0 = 2t$ and $\alpha^2 - \delta^2 = 4\beta$, 
the last expression becomes:
\begin{equation}\label{a^-_formula}
a^- = \frac{t}{2\delta}
  (\mathcal{T}_{\delta} - \mathcal{T}_{-\delta}t)
 +\frac{\alpha_0\delta^2 + \alpha\tau}{4\delta^2}
  (\mathcal{T}_\delta + \mathcal{T}_{-\delta} - 2I)
 +\frac{\beta}{\delta^2}X
  (\mathcal{T}_\delta + \mathcal{T}_{-\delta} - 2I). 
\end{equation}
Similarly, we have
\begin{align}
a^+ 
&= V - \frac{1}{2}a^0 
\nonumber\\
&= X - \frac{\alpha_0}{2}I 
 +\frac{1}{4}(X - \alpha_0)
   (\mathcal{T}_{\delta} + \mathcal{T}_{-\delta} - 2I)
 -\frac{1}{4\delta}(\alpha X + \tau)
   (\mathcal{T}_\delta - \mathcal{T}_{-\delta})
\nonumber\\
&\qquad -\frac{\alpha}{4\delta}(X - \alpha_0)
   (\mathcal{T}_\delta - \mathcal{T}_{-\delta})
 +\frac{\alpha}{4\delta^2}(\alpha X + \tau)
   (\mathcal{T}_\delta + \mathcal{T}_{-\delta} - 2I)
 -\frac{\alpha_0}{2}I 
\nonumber\\
&= X - \frac{\alpha}{2\delta}X
  (\mathcal{T}_{\delta} - \mathcal{T}_{-\delta}) 
 +\frac{\alpha^2 + \delta^2}{4\delta^2}X
  (\mathcal{T}_{\delta} + \mathcal{T}_{-\delta} - I)
\nonumber\\
&\qquad +\frac{\alpha\alpha_0 - \tau}{4\delta}
  (\mathcal{T}_\delta - \mathcal{T}_{-\delta}) 
 +\frac{\alpha\tau - \alpha_0\delta^2}{4\delta^2}
  (\mathcal{T}_{\delta} + \mathcal{T}_{-\delta} - I) 
 -\alpha_0I.
\label{V_2_Meixner_formula(2)}
\end{align}

We now introduce the {\em $\delta$-momentum operator} defined by
\begin{equation}
D_{\delta} = 
\begin{cases}
(\mathcal{T}_{\delta} - I)/\delta, & \text{if $\delta \neq 0$},\\
D & \text{if $\delta= 0$}.
\end{cases}
\end{equation}
%Let us observe that, due to the fact that $\mathcal{T}_c\mathcal{T}_d = \mathcal{T}_{c + d}$, for all $c$ and %$d$ complex numbers,
%we have:
%\begin{eqnarray}
%{\mathcal A}_{\delta/2}^2 & = & \frac{1}{\delta^2}\left(\mathcal{T}_{\delta/2} - {\mathcal %T}_{-\delta/2}\right)^2 \nonumber\\
%& = & \frac{1}{\delta^2}\left(\mathcal{T}_{\delta} + \mathcal{T}_{-\delta} - 2I\right),
%\end{eqnarray}
%if $\delta \neq 0$, and
%\begin{eqnarray}
%{\mathcal A}_{\delta/2}^2 & = & D^2,
%\end{eqnarray}
%if $\delta = 0$.\\
We have
\begin{equation}
D_{\delta}(F_n) \subseteq F_{n-1},
\qquad n \in {\mathbb N} \cup \{0\},
\end{equation}
which means that $D_{\delta}$ is a $(-1)$-faithful operator
with respect to the gradation $\{F_n\}_{n \geq 0}$.
Moreover, $D_{\delta}$ is a {\em purely momentuous operator}
in the sense that its position-momentum decomposition
contains no non-constant left-position coefficients $A_n(X)$.
Indeed,
\begin{equation}
D_{\delta} = \sum_{n=1}^{\infty}\frac{\delta^{n-1}}{n!}D^{n}.
\end{equation}
With this notation the annihilation operator of $X$ can be written as
\begin{equation}\label{05eqn:a- 5.60}
a^- 
= \frac{t}{2}(D_{\delta} + D_{-\delta}) 
 +\frac{\alpha_0\delta^2 + \alpha\tau}{4\delta}
 (D_{\delta} - D_{-\delta}) 
 +\frac{\beta}{\delta}X(D_{\delta} - D_{-\delta}).
\end{equation}
If $\beta = 0$, \eqref{05eqn:a- 5.60} becomes
\begin{equation}
a^- 
= \frac{t}{2}(D_{\delta} + D_{-\delta}) 
 +\frac{\alpha_0\delta^2 + \alpha\tau}{4\delta}
  (D_{\delta} - D_{-\delta}),
\end{equation}
which means that $a^-$ is a purely momentuous operator and 
is a linear combination of translation operators
whose coefficients add up to zero.
The above observation motivates us to study the following 
\begin{problem}\label{05problem:5.2}
Find all random variables $X$ having finite moments of all orders
such that its annihilation operator is of the form:
\begin{equation}\label{a^-_translation}
a^- = \sum_{i=1}^nc_i\mathcal{T}_{d_i}, 
\end{equation}
for some real constants $c_i$ and $d_i$, $1 \leq i \leq s$, 
satisfying the necessary condition:
\begin{equation}\label{cond_sum_coeff_zero}
\sum_{i = 1}^nc_i = 0. 
\end{equation}
Without loss of generality we may assume 
that $c_1, c_2, \dots, c_n$ are all different from zero
and $d_1, d_2, \dots, d_n$ are different each other.
\end{problem}

Before solving this problem, let us explain 
why the condition \eqref{cond_sum_coeff_zero} is necessary.
Since $a^-$ decreases the degree of a polynomial random variable
by $1$-unit, the degree of $a^-X^m$ must be at most $m-1$
for all $m \geq 0$.
Therefore the coefficient of $X^m$ in 
\begin{equation}
\sum_{i=1}^nc_i\mathcal{T}_{d_i}X^m 
=\sum_{i = 1}^nc_i(X + d_i)^m
\end{equation}
must be zero. 
This is equivalent to condition \eqref{cond_sum_coeff_zero}.

Going back to Problem \ref{05problem:5.2},
let us assume that $a^-$ has the form \eqref{a^-_translation}. 
We may assume that it has the mean $0$ after suitably shifting $X$
since a shift does not affect the annihilation operator $a^-$
Keeping in mind that $a^01 = E{X}1 = 0$, 
for all $m \geq 1$ we have
\begin{align}
E[X^m] 
&= \left\langle X^m1, 1\right\rangle 
\nonumber\\
&= \left\langle(a^+ + a^0 + a^-)X^{m-1}1, 1\right\rangle 
\nonumber\\
&= \left\langle X^{m-1}1, a^-1\right\rangle 
 + \left\langle X^{m-1}1, a^01\right\rangle
 + \left\langle a^-X^{m-1}1, 1\right\rangle 
\nonumber\\
&= 0 + 0 + 
 \left\langle \sum_{i=1}^nc_i\mathcal{T}_{d_i}X^{m-1}1,1\right\rangle
\nonumber\\
&= \sum_{i=1}^nc_i\left\langle(X+d_i)^{m-1}1,1\right\rangle
\nonumber\\
&= \sum_{i=1}^nc_iE[(X + d_i)^{m-1}].
\label{moment_recursion}
\end{align}

\noindent
{\bfseries Claim 1.} 
There exists a constant $k \geq 1$ such that 
\begin{equation}\label{moment_estimate_1}
|E[X^m]| \leq k^m\cdot m!,
\qquad m \geq 0.
\end{equation}
Indeed, from the obvious limit
\begin{equation}
\lim_{p \to \infty}\frac{|d_i|^p}{p!} = 0,
\qquad 1 \leq i \leq n,
\end{equation}
we see that there exists a constant $A \geq 1$ such that 
\begin{equation}
\frac{|d_i|^p}{p!} \leq A,
\qquad 1 \leq i \leq n, 
\qquad p \geq 0.
\end{equation}
We set
\begin{equation}
k = \max\left\{A\sum_{i = 1}^n|c_i|, 1\right\}.
\end{equation}
We now prove \eqref{moment_estimate_1} by induction on $n$.
Since $k \geq 1$, inequality \eqref{moment_estimate_1} is true
for $m = 0$.
Let $M\ge1$ and assume that inequality 
\eqref{moment_estimate_1} is true for all $m \leq M - 1$.
Using the recursive relation \eqref{moment_recursion}, we obtain
\begin{align}
|E[X^M]|
&= \left|\sum_{i=1}^nc_iE[(X + d_i)^{M-1}]\right|
\nonumber\\
&\leq \sum_{i=1}^n\sum_{m=0}^{M-1}|c_i|
 \binom{M-1}{m}|d_i|^{M-1-m}|E[X^m]|
\nonumber\\
&= \sum_{i=1}^n\sum_{m=0}^{M-1}|c_i|\binom{M - 1}{m}
 \cdot\frac{|d_i|^{M-1-m}}{(M-1-m)!}
 \cdot|E[X^{m}]| \cdot(M-1-m)!\,. 
\nonumber
\end{align}
Using the induction hypothesis, we obtain
\begin{align}
|E[X^M]|
&\leq \sum_{i=1}^n\sum_{m=0}^{M-1}|c_i|
 \binom{M-1}{m} \cdot A \cdot k^m m! \cdot (M-1-m)! 
\nonumber\\
&= \sum_{i=1}^n\sum_{m=0}^{M-1}|c_i|(M-1)!Ak^{m}. 
\nonumber
\end{align}
Note from $k \geq 1$ that $k^m \leq k^{M-1}$ for
for all $1\le m \le M-1$.
Then we have
\begin{align}
|E[X^M]|
&\leq \sum_{i=1}^n\sum_{m=0}^{M-1}|c_i|(M-1)!Ak^{M-1} 
\nonumber\\
&= \sum_{i=1}^nM|c_i|(M-1)!Ak^{M-1} 
\nonumber\\
&= k^{M-1}M!\left(A\sum_{i=1}^n|c_i|\right) 
\nonumber\\
&\leq k^MM!\,,
\end{align}
as desired.

\medskip
\noindent
{\bfseries Claim 2.} 
For all $m \in {\mathbb N} \cup \{0\}$ we have
\begin{equation}\label{moment_of_modulus_estimate}
E[|X|^m] \leq (2k)^mm!. 
\end{equation}
Indeed, since $X$ is a real-valued random variable, 
we have $|X|^2 = X^2$.
Using Jensen inequality for the concave function
$f(x) = \sqrt{x}$ ($x \geq 0$) and
the inequality \eqref{moment_estimate_1},
we have
\begin{align*}
E[|X|^m] 
&= E[\sqrt{X^{2m}}]
\leq \sqrt{E[X^{2m}]}
\leq \sqrt{k^{2m}(2m)!} \\
&=k^mm!\sqrt{\binom{2m}{m}}
\leq k^mm!\sqrt{(1 + 1)^{2m}}
= (2k)^mm!,
\end{align*}
as desired.

Multiplying both sides of \eqref{moment_recursion} by $t^{m-1}/(m-1)!$
and summing up for $m=1$ to $\infty$,
we obtain
\begin{equation}
\sum_{m=1}^{\infty}\frac{t^{m-1}}{(m-1)!}E[X^m]
=\sum_{i=1}^nc_i \sum_{m=1}^{\infty}
 \frac{t^{m-1}}{(m-1)!}E[(X+d_i)^{m-1}].
\end{equation}
We know that the above power series converges 
for $t$ in a neighborhood of $0$.
By the Lebesgue dominated convergence theorem
with \eqref{moment_of_modulus_estimate}
we can interchange the series and expectation
to obtain
\begin{equation}
E[X\exp(tX)]
=\sum_{i=1}^nc_iE[\exp(t(X + d_i))],
\end{equation}
where $t$ is in a neighborhood of $0$.
Let $\varphi(t)= E[\exp(tX)]$ be the Laplace transform of $X$.
Applying the Lebesgue dominated convergence theorem again,
we conclude that
\begin{equation}\label{diff_eq}
\frac{d\varphi}{dt}
=\varphi(t) \cdot \sum_{i = 1}^nc_i\exp(d_it). 
\end{equation}

We have two cases.

\noindent
{\bfseries Case 1.} $d_i \neq 0$ for all $1 \leq i \leq n$.
Since $\varphi(0) = E[\exp(0X)] = 1$, 
the solution of the differential equation \eqref{diff_eq} is
given by
\begin{equation}
\varphi(t) 
=\exp\left(\sum_{i=1}^n\frac{c_i}{d_i}\left(e^{d_it}-1\right)\right) 
=\exp\left(\sum_{i=1}^n\frac{c_i}{d_i}
  \left(e^{d_it}-d_it-1\right)\right),
\end{equation}
where $\sum_{i=1}^nc_i = 0$ is taken into account.
Then we have
\begin{equation}\label{Laplace_transform_product}
\varphi(t) 
=\prod_{i=1}^n
 \exp\left(\frac{c_i}{d_i}\left(e^{d_it}-d_it-1\right)\right). 
\end{equation}
On the other hand, 
letting $Y$ be a Poisson random variable with mean $\lambda>0$
and $\gamma$ a real constant,
we consider $Z= \gamma(Y - \lambda)$.
Then, for all $t \in\mathbb{R}$ we have
\begin{align}
E[\exp(tZ)] 
&= \sum_{k=0}^{\infty}
 e^{\gamma(k-\lambda)t}\frac{\lambda^k}{k!}e^{-\lambda} 
\nonumber\\
&= e^{-\lambda\gamma t - \lambda}\exp(\lambda e^{\gamma t})
\nonumber\\
&= \exp\left(\lambda(e^{\gamma t} - \gamma t - 1)\right).
\end{align}
Since according to formula \eqref{Laplace_transform_product}
the Laplace transform $\varphi(t)$ can be written
as a product of those of the above shifted 
and re-scaled Poisson random variables,
we conclude that there exist $n$ independent
Poisson random variables $Y_1, Y_2, \dots, Y_n$ with means
$c_1/d_1$, $c_2/d_2$, $\dots$, $c_n/d_n$, such that
\begin{equation}
X = \sum_{i = 1}^nd_iY_i.
\end{equation}
Since the mean (parameter) $\lambda$ of
a Poisson random variable must be positive, 
we have necessarily $c_i/d_i>0$ for all $1 \leq i \leq n$.
Hence $c_i$ and $d_i$ have the same sign for each $1 \leq i \leq n$.

\medskip
\noindent
{\bfseries Case 2.} $d_{i_0} = 0$ happens for some
$i_0 \in \{1,2,\dots, n\}$.
The solution of the differential equation \eqref{diff_eq} 
with initial condition $\varphi(0)=1$ is given by
\begin{align}
\varphi(t) 
&= \exp\left(c_{i_0}t+\sum_{i\neq i_0}\frac{c_i}{d_i}
  \left(e^{d_it} - 1\right)\right)
\nonumber\\
&= \exp\left(\sum_{i \neq i_0}\frac{c_i}{d_i}
  \left(e^{d_it}-d_it-1\right)\right)
\nonumber\\
&= \prod_{i \neq i_0}
   \exp\left(\frac{c_i}{d_i}\left(e^{d_it}-d_it-1\right)\right), \label{new_Laplace_transform_product}
\end{align}
where $\sum_{i \neq i_0}c_i = -c_{i_0}$ is taken into account.
Thus
\begin{equation}
X=\sum_{i \neq i_0}d_iY_i\,,
\end{equation}
where $Y_i$, $1\le i\le n$, $i\neq i_0$,
are independent Poisson random variables with means $c_i/d_i$.

Since we have initially shifted the random variable $X$ to have
the mean equal to zero, we can say that the answer to the problem
of describing the random variables having finite moments of all orders,
whose annihilation operator is a linear combination of translation
operators, is that those random variables are all shifted linear combinations of independent Poisson random variables.

Upon closing this paper we note that 
formula \eqref{a^-_formula} can be used 
to recover all the six types of Meixner random variables, 
by centering first the random variable $X$ 
and then using the formula:
\begin{align}
E[X^m] 
&= \langle X^m1,1\rangle 
\nonumber\\
&= \langle (a^+ + a^0 + a^-)X^{m-1}1, 1\rangle
\nonumber\\
&= \langle X^{m-1}1, a^-1 \rangle 
 +\langle X^{m-1}1, a^01 \rangle
 +\langle a^-X^{m-1}1,1 \rangle
\nonumber\\
&= 0 + 0 + \langle a^-X^{m-1}1, 1\rangle
\nonumber\\
&= \langle a^-X^{m-1}1,1\rangle.
\end{align}
Substituting $a^-$ by its position-momentum decomposition 
\eqref{a^-_formula}, 
we obtain first a recursive formula for the moments of $X$ as before,
then a differential equation for the Laplace transform of $X$.
Solving this differential equation and inverting the Laplace transform,
we obtain the well-known classification of
the classical Meixner random variables as below.

\begin{enumerate}
\item $\alpha = \beta = 0$.
In this case $X$ is a {\em Gaussian} random variable, 
i.e., a continuous random variable with the density function
\begin{equation*}
f_X(x) = \frac{1}{\sqrt{2\pi t}}\,
 \exp\bigg(-\frac{(x - \alpha_0)^2}{2t}\bigg).
\end{equation*}
\item $\beta = 0$ and $\alpha \neq 0$.
In this case $X$ is a shifted and re-scaled {\em Poisson}
random variable and its distribution is given by
\begin{equation*}
\mu_X 
= \sum_{k=0}^{\infty}
 \frac{\lambda^k}{k!}\,
 e^{-\lambda}\delta_{\alpha(k -\lambda) + \alpha_0},
\end{equation*}
where $\lambda=t/\alpha^2$.
\item $\beta > 0$ and $\alpha^2 > 4\beta$.
Then $X$ is a shifted {\em Pascal (negative binomial)}
random variable with the distribution given by
\begin{equation*}
\mu_X 
=\sum_{k=0}^{\infty} \frac{\Gamma(r+k)}{k!\Gamma(r)}\,
 p^r(1-p)^k\delta_{k-[2t/(\alpha + d)]+\alpha_0},
\end{equation*}
where
\[
d=\sqrt{\alpha^2 - 4\beta},
\qquad p = \frac{2d}{\alpha + d}\,,
\qquad
r=\frac{t}{\beta}\,.
\]
\item $\beta > 0$ and $\alpha^2 = 4\beta$.
In this case $X$ is a shifted and re-scaled 
{\em Gamma} random variable with shift parameter
$2t/\alpha$ and scaling parameter $\alpha/2$, i.e.,
its density function is given by
\begin{equation*}
f_X(x)
=\frac{2^{2t/\alpha}}{\alpha^{2t/\alpha}\,\Gamma(2t/\alpha)}\,
x^{(2t/\alpha)- 1}e^{-2x/\alpha}1_{(0, \infty)}(x).
\end{equation*}
\item $\beta > 0$ and $\alpha^2 < 4\beta$.
In this case,
up to translation $X$ is a {\em two parameter hyperbolic secant} 
random variable whose density function is given by
\begin{equation*}
f_X(x) = ce^{2\theta x/\gamma}\left|\Gamma(k + ix\gamma)\right|^2,
\end{equation*}
where
\[
\gamma=\sqrt{4\beta - \alpha^2}\,,
\qquad
k = \frac{2t}{r\gamma}\,,
\qquad
\gamma + i\alpha = re^{i\theta},
\quad -\frac{\pi}{2} < \theta < \frac{\pi}{2}\,.
\]
\item $\beta < 0$. Then $t \in -{\mathbb N}\beta$.
In this case, up to a shifting and re-scaling 
$X$ is a {\em binomial} random variable
whose distribution is given by
\begin{equation*}
\mu_X =\sum_{k = 0}^n\binom{n}{k}p^k(1 - p)^{n - k}\delta_k,
\end{equation*}
where 
\[
n = -\frac{t}{\beta}\,,
\qquad
p = \frac{1 \pm \sqrt{c/(4 + c)}}{2}\,,
\qquad
c=-\frac{\alpha^2}{\beta} \geq 0.
\]
\end{enumerate}

{\bfseries Acknowledgements.} 
This research started during the visits of Aurel I. Stan 
to Ochanomizu University in Tokyo and Tohoku University in Sendai, 
Japan in November 2019,
and developed during the stay of Nobuaki Obata
at Ohio State University in June 2023.
We would like to thank the JSPS Grant-in-Aid for Scientific Research
XXXXXX and 19H01789 for supporting these research visits.

\end{document}